\documentclass[reqno, 11pt]{amsart}
\usepackage{amsmath,amssymb,amsthm}

\usepackage[
    a4paper, 
    hmargin={2.7cm,2.7cm}, 
    vmargin={3.3cm,3.3cm}
]{geometry}
\usepackage[colorlinks=true,linkcolor=blue,citecolor=blue]{hyperref}

\usepackage[
    backend    = biber,
    sortcites  = false, 
    giveninits = true,
    doi        = false,
    isbn       = false,
    url        = false,
    maxnames   = 50,
    citestyle  = alphabetic,
    bibstyle   = alphabetic
]{biblatex}
\DeclareNameAlias{sortname}{family-given}
\DeclareNameAlias{default}{family-given}
\theoremstyle{plain}
\newtheorem{maintheorem}{Theorem}

\newtheorem{theorem}{Theorem}[section]
\newtheorem{corollary}[theorem]{Corollary}
\newtheorem{lemma}[theorem]{Lemma}
\newtheorem{proposition}[theorem]{Proposition}
\theoremstyle{definition}

\newtheorem{remark}[theorem]{Remark}
\newtheorem{example}[theorem]{Example}

\newcommand{\Q}{\mathbb{Q}}
\newcommand{\R}{\mathbb{R}}
\newcommand{\C}{\mathbb{C}}
\newcommand{\dif}{\,\mathrm{d}}
\newcommand{\cG}{\mathcal{G}}
\newcommand{\cP}{\mathcal{P}}
\newcommand{\SL}{\mathrm{SL}}
\newcommand{\m}[1]{ \begin{pmatrix} #1 \end{pmatrix} }

\DeclareMathOperator{\supp}{supp}

\date{November 18, 2026}
\title
[Boundaries of compact-by-amenable groups]
{Boundaries of groups with compact-by-amenable decomposition}

\author{Hannes Wendt}
\email{hhwendt@math.uio.no}
\address{Department of Mathematics, University of Oslo}

\author{Makoto Yamashita}
\email{makotoy@math.uio.no}
\address{Department of Mathematics, University of Oslo}

\makeatletter
\def\@setaddresses{} 
\let\@\@maketitle\def\@maketitle{\@ 
  \def\author##1##2{\par##2}
  \def\email##1##2{Email: \texttt{##2}, }
  \def\address##1##2{Address: ##2 }
  \@footnotetext{\addresses} 
}
\makeatother

\begin{document}
\renewcommand{\epsilon}{\varepsilon}

\begin{abstract}
We look at the Liouville property (triviality of harmonic functions) of topological groups $G$ relative to a compact subgroup $K$.
We show that Gelfand pairs have this property.
Furthermore, we show that this property holds if and only if $G$ admits a product decomposition $G = K S$ for an amenable subgroup $S$.
We identify the Poisson and Furstenberg boundaries of these groups.
This applies in particular to reductive algebraic groups over local fields.
\end{abstract}

\maketitle\thispagestyle{empty}

\section{Introduction}

Let $\mu$ be a probability measure on a locally compact group $G$, and let $K \leq G$ be a compact subgroup. We say \emph{$\mu$ has the Liouville property relative to $K$} if the following condition holds.
\begin{quote}
For any bounded left-$K$-invariant function $f$ on $G$, the equality $f = \int_G f(\,\cdot\, h) \dif \mu(h)$ implies that $f$ is constant.
\end{quote}
Similarly, we say that $(G, K)$ has the Liouville property if there is such an admissible probability measure $\mu$.
We show that this property holds whenever $(G, K)$ is a Gelfand pair, in the following sense.

\begin{maintheorem}\label{thmA:Liouville-from-comm-Hecke-alg}
Let $G$ be a locally compact second countable group, and let $K \leq G$ be a compact subgroup such that the associated Hecke algebra $C_c(K \backslash G / K)$ is commutative. Then all admissible, bi-$K$-invariant probability measures have the Liouville property relative to $K$.
\end{maintheorem}

For $K = \{e\}$ the above recovers the Choquet--Deny theorem~\cite{ChDe60}.
Our direct motivation comes from Furstenberg's celebrated work~\cite{Fu63} on the boundaries of semisimple Lie groups $G$, in which a result of the above type for the maximal compact subgroup $K \leq G$ plays a crucial role in the identification of the Poisson boundary of $G$ through analysis of random walks on the symmetric space $G / K$.

Indeed, the relative Liouville property allows us to get a strong structural property for locally compact groups and their Poisson boundaries closely analogous to semisimple Lie groups.
More precisely, we have the following results.

\begin{maintheorem}\label{thmB:Compact-times-amenable}
    Let $G$ be a locally compact second countable group, and let $K \leq G$ be a compact subgroup. \begin{enumerate}
        \item If $G = KS$ holds for some closed amenable subgroup $S \leq G$, then there exists an admissible probability measure on $G$ that has the Liouville property relative to $K$.
        \item If there exists an admissible probability measure $\mu$ on $G$ that has the Liouville property relative to $K$, then $G = KS$ for a closed, amenable subgroup $S \leq G$, and the Poisson boundary of $(G,\mu)$ is $(G/S, \nu)$ for some $G$-quasi-invariant probability measure $\nu$.
    \end{enumerate}
\end{maintheorem}

For $K = \{e\}$ the above recovers the main result of~\cite{Ros81}, which states that amenability of $G$ is characterized by the triviality of some Poisson boundary of $G$.
More recently, Monod proved that Gelfand pairs $(G,K)$ admit an Iwasawa decomposition~\cite{Mo20}. The combination of Theorems~\ref{thmA:Liouville-from-comm-Hecke-alg} and~\ref{thmB:Compact-times-amenable} gives a random-walk proof of his result. 

Finally, we also have a description of the Furstenberg boundary for groups $G$ satisfying the conclusion of the above theorem, as follows.

\begin{maintheorem}\label{thmC:Homogeneous-Poisson-is-Furstenberg}
    Let $G$ be a locally compact second countable group such that some Poisson boundary of $G$ is isomorphic to $(G/S, \nu)$ where $S$ is a cocompact, amenable subgroup and $\nu$ a $G$-quasi-invariant probability measure. Then $G/S$ is the Furstenberg boundary of $G$.
\end{maintheorem}

A prominent class of locally compact groups satisfying the assumption of Theorem~\ref{thmA:Liouville-from-comm-Hecke-alg} comes from connected reductive algebraic groups $\cG$ over non-Archimedean local fields $F$, by taking the group of $F$-points $G = \cG(F)$.
Using Theorems~\ref{thmB:Compact-times-amenable} and~\ref{thmC:Homogeneous-Poisson-is-Furstenberg}, we identify the Furstenberg boundary of such groups $G$, see Section~\ref{sec:red-alg-grps}.

As a concrete example, we can consider the determinant-one matrices over the $p$-adic numbers $G = \SL_n(\Q_p)$ and identify its Furstenberg boundary as $G/B$, where $B$ is the Borel subgroup of upper triangular matrices.
This is a direct analogue of the well-known result in the context of semisimple Lie groups established in Furstenberg's foundational work~\cite{Fu73}.

While this paper was in the final stage of preparation, Monod announced~\cite{Mon26}. His work contains closely overlapping results particularly on the Furstenberg boundary of groups $G$ as above, although the proof strategy seems somewhat different from ours.

\bigskip
\paragraph{\bf Acknowledgements}
Part of this work was carried out during the first-named author's stay at Institut Henri Poincaré for the \emph{Operator algebras: approximation, rigidity and dynamics} program in the fall of 2026.
We thank Sven Raum for a stimulating discussion.

\noindent M.Y.'s research was funded, in part, by The Research Council of Norway [project 324944].

\smallskip
\paragraph{Declaration on AI usage}
After developing the main ideas for this paper, we used a large language model for part of the brainstorming process to refine formulations.

\section{Preliminaries}

Let $G$ be a locally compact, second countable group.
Throughout the paper we fix a left Haar measure on $G$, denoted $m_G$.

\subsection*{Poisson boundary}
Let $\mu$ be a probability measure on $G$.
We say that:
\begin{enumerate}
    \item $\mu$ is \emph{spread out} if $\mu^{*n} \ll m_G$ for some $n \geq 1$, that is, some positive convolution power of $\mu$ is absolutely continuous relative to the Haar measure on $G$.
    \item $\mu$ is \emph{adapted} if $\supp \mu$ generates $G$ as a semigroup, that is, $G = \bigcup_{n \geq 0} (\supp \mu)^n$.
    \item $\mu$ is \emph{admissible} if $\mu$ is spread out and adapted.
\end{enumerate} 

Let $X$ be a standard Borel $G$-space.
For a probability measure $\eta$ on $X$, our convention of the convolution is $\mu * \eta = \int_G \eta(g^{-1}\,\cdot\,) \dif \mu(g)$, which is a new probability measure on $X$.
We say $\eta$ is \emph{$\mu$-stationary} if $\eta = \mu * \eta$. By~\cite[Lemma 1.1]{NZ99}, if $\mu$ is admissible, then $\eta$ is $G$-quasi-invariant, that is $\eta(g \,\cdot\,)$ and $\eta$ have the same null sets on $X$ for all $g \in G$.

Therefore given an admissible probability measure $\mu$ on $G$, for every standard Borel $G$-space $(X, \eta)$ where $\eta$ is a $\mu$-stationary measure, we have a well-defined map between von Neumann algebras \[
    P_\eta\colon L^\infty(X, \eta) \to L^\infty(G, m_G), \quad P_\eta(F)(g) = \int_X F(gx) \dif \eta(x).
\]
Automatically, $P_\eta$ is linear, unital, positive, normal, contractive and $G$-equivariant. The image of $P_\eta$ can be as small as $\C 1_G$, and by evaluating indicator functions one finds that $P_\eta$ is isometric if and only if $\operatorname{ess\,sup}_{g \in G} \eta(g^{-1}A) = 1$ holds for all Borel set $A \subseteq X$ satisfying $\eta(A) > 0$.

By~\cite[Theorem 3.1]{Fu71}, there exists a standard Borel $G$-space $(\Pi_\mu, \hat\mu)$, called the \emph{Poisson boundary of $(G,\mu)$}, which has the property that $P_{\hat\mu}$ is isometric and the image of $P_{\hat\mu}$ is the set of the \emph{$\mu$-harmonic} functions $\{f \in L^\infty(G, m_G) \colon f = f * \mu\}$. Here $f * \mu = \int_G f(\,\cdot\, h) \dif \mu(h)$ is an established but unfortunate notation, because $f * (\mu_1 * \mu_2) = (f * \mu_2) * \mu_1$ reverses the order of convolution, where $\mu_1, \mu_2$ are admissible measures on $G$.

There is a separable, metrizable topological $G$-space that recovers $\Pi_\mu$ by taking its Borel $\sigma$-algebra, see for example~\cite[Theorem 2.1.19]{Zi84}.

\subsection*{Furstenberg boundary}

Let $X$ be a compact space with a continuous $G$-action. We say that:
\begin{enumerate}
    \item $X$ is \emph{minimal} if it has no closed $G$-invariant subspaces, except $\emptyset$ and $X$.
    \item $X$ is \emph{strongly proximal} if for every probability measure $\eta$ on $X$ there is a point mass contained in the weak$^*$-closure of the $G$-orbit $G\eta$.
    \item $X$ is \emph{a $G$-boundary} if it is minimal and strongly proximal.
\end{enumerate}
By~\cite[Proposition 4.6]{Fu73} there exists a \emph{Furstenberg boundary} or \emph{universal boundary} $\partial_F G$, that is an up to isomorphism unique $G$-boundary with the universal property that for any $G$-boundary $X$ there exists a $G$-equivariant surjection $\partial_F G \twoheadrightarrow X$.

\section{Liouville measures from Gelfand pairs}

In this section we prove Theorem~\ref{thmA:Liouville-from-comm-Hecke-alg}. 
We recall that $(G,K)$ is a Gelfand pair, that means the Hecke algebra $C_c(K \backslash G / K)$ is commutative, if and only if all bi-$K$-invariant measures on $G$ commute. For any $a \in G$ the measure $m_{KaK} = m_K * \delta_a * m_K$ is an example of a bi-$K$-invariant measure, where $m_K$ denotes the Haar probability measure on the compact subgroup $K$.

The only observation we need to make is that whenever $\mu$ is a bi-$K$-invariant probability measure and $f$ is a $\mu$-harmonic function, also $f * m_{KaK}$ are $\mu$-harmonic for all $a \in G$: \[
    (f * m_{KaK}) * \mu
    = f * (\mu * m_{KaK})
    = f * (m_{KaK} * \mu)
    = (f * \mu) * m_{KaK}
    = f * m_{KaK}.
\]
Also, for all bi-$K$-invariant $\nu$ on $G$ we may write \[
    f * \nu 
    = \int_G f(\,\cdot\, h) \dif \nu(h)
    = \int_G (f * \delta_h) \dif \nu(h)
    = \int_G (f * m_{KhK}) \dif \nu(h).
\]
 
\begin{proof}[Proof of Theorem~\ref{thmA:Liouville-from-comm-Hecke-alg}]
    Since $\mu$ is spread out, there exists $\rho \in L^1(G, m_G)$ for some $m \geq 1$ such that $\mu^{*m} = \rho \cdot m_G$. So any $\mu$-harmonic $f \in L^\infty(G, m_G)$ has a continuous representative \[
        f(g) = (f * \mu^{*m})(g) = \int_G f(g h) \rho(h) \dif m_G(h) = \int_G (f \cdot \lambda_g \rho) \dif m_G.
    \]
    Also, we recall that $\mu$-harmonic $f$ is automatically right-$K$-invariant because $\mu = \delta_k * \mu$.
    
    Consider \[
        \{f \mid \text{bounded, left-$K$-invariant, $\mu$-harmonic}\}
        = \operatorname{span} C
    \] where \[
        C = \left\{f \colon G \to [0,1] \mid 
        \forall k \in K \colon f = \lambda_k f \text{ and }
        f = f * \mu \right\}
        \subseteq L^\infty(G, m_G)
    \] is convex and bounded. Because $\lambda_k$ and $f \mapsto f * \mu$ are weak$^*$-continuous, the latter being the Banach-space adjoint of a contraction on $L^1(G, m_G)$, we find that $C$ is weak$^*$-closed. By the Banach--Alaoglu theorem $C$ is weak$^*$-compact. So by the Krein--Milman theorem, it suffices to show the extreme points of $C$ are constant.

    Consider an extreme point $f \in C$. For $\nu = \mu^{*n}$ and $t = \nu(A) \in (0,1)$ we have \[
        f = f * \nu = t \left( \frac1{\nu(A)} \int_A (f * m_{KaK}) \dif \nu(a) \right) + (1 - t) \left( \frac1{\nu(G \setminus A)} \int_{G \setminus A} (f * m_{KaK}) \dif \nu(a) \right)
    \] Because $f * m_{KaK} \in C$, and limits of convex combinations remain in $C$, extremality implies \[
        f = \frac1{\nu(A)} \int_A (f * m_{KaK}) \dif \nu(a), \quad \nu(A) > 0.
    \] Note that for $\nu(A) = 1$ this holds trivially because $f = f * \nu$.
    By continuity, we may evaluate both sides at $e$ and use the bi-$K$-invariance of $f$ in the form $(f * m_{KaK})(e) = f(a)$ to obtain that $f(e) = \frac1{\nu(A)} \int_A f(a) \dif \nu(a)$ for all $\nu(A) > 0$.
    This implies both $\{f > f(e)\}$ and $\{f < f(e)\}$ must be disjoint from $\supp \nu$. Since $\mu$ is adapted, this shows $f$ is constant on all of $G$.
\end{proof}

\begin{remark}
    In the case when $K$ is a compact, open subgroup of $G$, the countable set measures $m_{KaK}$ over double-cosets $KaK \in K \backslash G / K$ forms a hypergroup under convolution. Combining results of~\cite{Yam99} and~\cite{Hay00}, commutative hypergroups with dimension function can be shown to be weakly amenable in the sense of~\cite{HiIz98}, which is witnessed by an ergodic probability measure. By an adaptation of this proof, on the group $G$ one can show that in fact all absolutely continuous, bi-$K$-invariant probability measures on $G$ with full support have the Liouville property relative to $K$.
\end{remark}

As we remarked in the introduction, connected reductive algebraic groups over totally disconnected local fields satisfy the assumption for Theorem~\ref{thmA:Liouville-from-comm-Hecke-alg}, see Section~\ref{sec:red-alg-grps}.

\section{Liouville measure from Iwasawa decomposition}

In this section, we prove the part of Theorem~\ref{thmB:Compact-times-amenable} that constructs an admissible probability measure with the Liouville property from a compact-by-amenable decomposition. 
We closely follow the proof of~\cite[Theorem 4.3]{KaVe83}. The only additional ingredient is the following:

\begin{lemma}\label{lemma:reiter}
    Assume $G = KS$ with $K$ compact and $S$ closed and amenable. Then for every compact subset $C \subseteq G$ and $\epsilon > 0$ there exists a probability measure $\alpha$ on $G$ such that \[
        \|m_K * \delta_g * \alpha - m_K * \alpha\| < \epsilon, \quad g \in C.
    \]
\end{lemma}

Here $\|\cdot\|$ is always the total variation norm for measures on $G$.

\begin{proof}[Proof of Lemma~\ref{lemma:reiter}]
    The set $Q = S \cap KC$ is compact. By Reiter's property, using the amenability of $S$, there exists a probability density $\varphi$ on $S$ such that $\alpha = \varphi \, m_S$ satisfies $\|\delta_s * \alpha - \alpha\| = \int_S |\varphi(s^{-1} t) - \varphi(t)| \dif m_S(t) < \epsilon$ for all $s \in Q$. Here $m_S$ is the Haar measure on $S$. Now for all $g \in C$ we have a decomposition $g = ks$ where $k \in K$ and $s = k^{-1}g \in S \cap KC = Q$, and $m_K * \delta_g = m_K * \delta_s$ by the $K$-invariance of the Haar measure on $K$. Hence \[
        \|m_K * \delta_g * \alpha - m_K * \alpha\|
        = \|m_K * (\delta_s * \alpha - \alpha)\|
        \leq \|\delta_s * \alpha - \alpha\| < \epsilon.
        \qedhere
    \]
\end{proof}

To construct a Liouville measure relative to $K$ we only have to construct $\mu$ such that \begin{equation}\label{eq:sufficient-liouville}
    \lim_{n \to \infty} \|m_K * \delta_g * \mu^{*n} - m_K * \mu^{*n}\| = 0, \quad \text{ for all } g \in G.
\end{equation}
This is because for all left-$K$-invariant, $\mu$-harmonic $f$ we have $f(g) = (m_K * \delta_g * \mu^{*n})(f)$.
So for all bounded, left-$K$-invariant, $\mu$-harmonic functions $f$, we get an estimate \[
    |f(g) - f(e)| = |(m_K * \delta_g * \mu^{*n} - m_K * \mu^{*n})(f)| \leq \|f\|_\infty \|m_K * \delta_g * \mu^{*n} - m_K * \mu^{*n}\| \to 0
\]
for all $g \in G$, from~\eqref{eq:sufficient-liouville}.

There is a generic argument from the proof of~\cite[Theorem 4.3]{KaVe83} which shows how to combine probability measures $\alpha_1, \alpha_2, \dots$ obtained from Lemma~\ref{lemma:reiter} into a single probability measure \[
    \mu = p_0 \alpha_0 + p_1 \alpha_1 + \dots,
\] for certain weights $p_k > 0$ with $\sum_{k \geq 0} p_k = 1$ such that~\eqref{eq:sufficient-liouville} holds.
Choosing $\alpha_0$ with strictly positive density relative to $m_G$ gives $\mu$ full support, and all $\alpha$ from Lemma~\ref{lemma:reiter} can be chosen absolutely continuous by convolving on the right with some fixed probability $\mu_0 \ll m_G$.

\begin{proof}[Proof of part (1) in Theorem~\ref{thmB:Compact-times-amenable}]
    We pick compact sets $e \in C_1 \subseteq C_2 \subseteq \dots$ with $\bigcup_{k \geq 1} C_k = G$. Pick $p_k, \epsilon_k > 0$ with $\sum_{k \geq 0} p_k = 1$ and $\epsilon_k$ decreasing to $0$ and such that $(p_0 + \dots + p_{n-1})^n < \epsilon_n$ for all $n \geq 1$. Now, for any $k \geq 1$, given $\alpha_0, \dots, \alpha_{k-1}$ consider the finite set of probability measures \[
        F_k = \{ \alpha_{w_1} * \dots * \alpha_{w_m} \mid 0 \leq w_j < k \text{ for all } 0 \leq j \leq m < k\}.
    \] We may select (e.g. among the $C_k$) a large enough compact set $D_k$ such that for all $\beta \in F_k$ it holds $\beta(G \setminus D_k) < \epsilon_k$. Now finally, we pick using Lemma~\ref{lemma:reiter} a probability measure $\alpha_k$ such that \[
        \| m_K * \delta_h * \alpha_k - m_K * \alpha_k \| < \epsilon_k, \quad h \in (C_k D_k).
    \]
    Now for $\mu = \sum_{k \geq 0} p_k \alpha_k$ write $p_w = p_{w_1} \dots p_{w_n}$ and $\alpha_w = \alpha_{w_1} * \dots * \alpha_{w_n}$ for words $w \in \mathbb N^n$. Then \[
        \mu^{*n}
        = \sum_{w \in \mathbb N^n} p_w \alpha_w
        = \sum_{w \,\in\, \{0, \dots, n-1\}^n} p_w \alpha_w + \sum_{
            \substack{w \,\in\, \mathbb N^n \setminus \{0, \dots, n-1\}^n \\ \Rightarrow\, (\exists k \geq n \text{ and } \beta \in F_n)}
        } p_w (\beta * \alpha_k * \gamma).
    \] Here the first summand contains only words where all the $n$ letters are less than $n$. The second summand contains the rest of all words, where there is some letter $k$ that is at least $n$ and we choose the first such occurrence. So there $\alpha_w = \beta * \alpha_k * \gamma$ where $\beta \in F_n \subseteq F_k$ and $\gamma$ is some probability measure made of the remaining factors $\alpha_{w_j}$.
    For all $g \in G$ it holds \[\begin{aligned}
        &\| m_K * \delta_g * \mu^{*n} - m_K * \mu^{*n} \|
        \leq \sum_{w \in \mathbb N^{n}} p_w \| m_K * \delta_g * \alpha_w - m_K * \alpha_w \| \\
        &\quad\leq 2 \sum_{w \,\in\, \{0, \dots, n-1\}^n} p_w + \sum_{\text{the rest}} p_w \|m_K * \delta_g * \beta * \alpha_k - m_K * \beta * \alpha_k\| \|\gamma\| \\
        &\quad\leq 2 \epsilon_n
        + \sup_{k \geq n, \, \beta \in F_n} \int_G \|m_K * \delta_g * \delta_x * \alpha_k - m_K * \delta_x * \alpha_k \| \dif \beta(x)
    \end{aligned}\]
    Here is the only place we use the choice of $p_k$, to ensure the first sum over the $n^n$ words is small, where we have used the crude bound $2$ for the difference of probability measures.
    
    For the integral expression, we have
    \[
        \int_G \| m_K * \delta_{gx} * \alpha_k - m_K * \delta_{ex} * \alpha_k \| \dif \beta(x)
        \leq \int_{G \setminus D_k} 2 \dif \beta + \int_{D_k} 2\epsilon_k \dif \beta
        \leq 4\epsilon_k
    \]
    since the case $x \in D_k$ allows us to replace $m_K * \delta_{gx} * \alpha_k$ by $m_K * \alpha_k$ with an $\epsilon_k$-error, while the case $x \not\in D_k$ has small $\beta$-measure bounded by $\epsilon_k$ again.
    So we have shown $\| m_K * \delta_g * \mu^{*n} - m_K * \mu^{*n} \| < 6 \epsilon_n$ for every $g \in C_n$.
\end{proof}

\section{Poisson boundary from a Liouville measure}

We are now ready to give applications to boundary theory.
In this section, we prove the part of Theorem~\ref{thmB:Compact-times-amenable} that constructs the Iwasawa decomposition and Poisson boundary.

\begin{proof}[Proof of part (2) in Theorem~\ref{thmB:Compact-times-amenable}]
    Let $P = \supp \hat\mu \subseteq \Pi_\mu$. It is closed, and $G$-invariant because $\hat\mu$ was $G$-quasi-invariant. Because $\Pi_\mu \setminus P$ is a null set, we have a $G$-equivariant isomorphism of measure spaces \[
        (\Pi_\mu, \hat\mu) \cong (P, \hat\mu|_P).
    \]

    Similarly as in~\cite[Lemma 5.1]{Fu63},
    for every bounded continuous $\phi \in C_b(\Pi_\mu)$, we consider the function \[
        \phi'(\xi) = \int_K \phi(k\xi) \dif m_K(k), \qquad \xi \in \Pi_\mu.
    \]
    It is bounded, and $K$-invariant by right invariance of $m_K$ (compact groups are unimodular). By continuity of the action and dominated convergence, $\phi'$ is continuous.
    Then $P_{\hat\mu}(\phi')$ is a bounded, left-$K$-invariant $\mu$-harmonic function. By the Liouville property of $\mu$ it is constant, say $P_{\hat\mu}(\phi') = c \in \C$. Because $P_{\hat\mu}$ is linear, unital and injective, $\phi' = c$ holds $\hat\mu$-a.e. So the set $\{\phi' \neq c\}$ is an open null set. By definition, it must be disjoint from $P$. This shows \begin{equation}\label{eq:phi-prime-constant}
        \forall \phi \in C_b(\Pi_\mu)\colon \phi' \text{ is constant on } P.
    \end{equation}

    Fix $x \in P$. By $G$-invariance of $P$ we have $Kx \subseteq P$. Suppose there is $y \in P$ with $Ky \neq Kx$. Then the $K$-orbits $Kx$ and $Ky$ are disjoint compact sets, and we find $\phi \in C_b(\Pi_\mu)$ with $\phi|_{Kx} = 0$ and $\phi|_{Ky} = 1$. Because the topology on $\Pi_\mu$ is metrizable, one does not even need Urysohn's lemma but can construct $\phi$ concretely. In any case $\phi'(x) = 0$ and $\phi'(y) = 1$, contradicts~\eqref{eq:phi-prime-constant}, so $P = Kx$. Since $P$ is $G$-invariant, we have $Gx \subseteq P = Kx \subseteq Gx$, so we have shown \[
        P = Gx = Kx.
    \]
    Define \[
        S = G_{x} = \{g \in G \mid gx = x\}.
    \] For all $g \in G$ we know that $gx \in Kx$, so $g x = k_g x$ for some $k_g \in K$, that is $k_g^{-1} g \in S$, i.e.\ $g \in KS$. That means $G = KS$. Importantly, the compact $G/S$ is Hausdorff because $S$ is closed.
    Now the orbit map \[
        \psi \colon G/S \to P, \quad gS \mapsto gx,
    \] is a well-defined $G$-equivariant bijection. It is continuous by the universal property of the quotient topology, and it is actually a homeomorphism because continuous bijections from compact spaces onto Hausdorff spaces are closed.
    So we get a $G$-equivariant isomorphism of measure spaces \[
        (\Pi_\mu, \hat\mu) \cong (P, \hat\mu|_P) 
        \stackrel{\psi^{-1}}\longrightarrow (G/S, \nu), 
        \quad \nu = (\psi^{-1})_*(\hat\mu|_P).
    \]
    By equivariance, again $\nu$ is $G$-quasi-invariant.

    The action $G \curvearrowright (\Pi_\mu, \hat\mu)$ is Zimmer-amenable~\cite[Theorem 5.2]{Zi78}. The Zimmer-amenability does not feel null sets or equivariant isomorphisms of measure spaces, so $G \curvearrowright (G/S, \nu)$ is Zimmer-amenable. By~\cite[Proposition 4.3.2]{Zi84} this is equivalent to $S$ being amenable.
\end{proof}

The following is essentially the computation leading up to~\cite[Theorem 5.4]{Fu63}.

\begin{corollary}
    In the setting of Theorem~\ref{thmB:Compact-times-amenable} part (2), if in addition $\mu$ is left-$K$-invariant, then $\nu$ is the unique $K$-invariant probability measure on $G/S$. Namely, $\nu = m_K * \delta_{eS}$.
\end{corollary}
\begin{proof}
    If $\mu$ is left-$K$-invariant, then $\hat\mu$ is $K$-invariant, so also $\nu$ is $K$-invariant. 
    Let $\nu'$ be any $K$-invariant probability measure on $G/S$ and $f \in C(G/S)$. By $K$-invariance and Fubini \[
        \nu'(f) 
        = \int_K \int_{G/S} f(kz) \dif\nu'(z) \dif m_K(k)
        = \int_{G/S} \int_K f(kz) \dif m_K(k) \dif\nu'(z).
    \] Every $z \in G/S$ is of the form $z = k_0 S$ with $k_0 \in K$, because $G = KS$. By right-$K$-invariance of $m_K$ (compact groups are unimodular) we have $\int_K f(kk_0S) \dif m_K(k) = \int_K f(kS) \dif m_K(k)$. Hence $\nu'(f) = \int_K f(kS) \dif m_K(k)$.
\end{proof}

\section{Furstenberg boundary from homogeneous Poisson boundary}

Let us next prove Theorem~\ref{thmC:Homogeneous-Poisson-is-Furstenberg}, which gives an application to the Furstenberg boundary.

From~\cite[Proposition 6.4]{CKM23} we know that $G$ has a cocompact amenable subgroup if and only if $G$ acts transitively on the Furstenberg boundary. By~\cite[Proposition 10 and the preceding discussion]{Oza14}, in the setting of Theorem~\ref{thmC:Homogeneous-Poisson-is-Furstenberg} we know already that the Furstenberg boundary of $G$ must be of the form $G/\widetilde S$ for some maximal cocompact, relatively amenable subgroup $\widetilde S$ that is unique up to conjugation. 

In order to prove Theorem~\ref{thmC:Homogeneous-Poisson-is-Furstenberg}, namely that $\widetilde S$ can be taken to be $S$, we have to show:

\begin{proposition}\label{prop:G-boundary-from-poisson}
    Suppose some Poisson boundary of $G$ is isomorphic to $(G/H, \nu)$ with some closed, cocompact subgroup $H \leq G$ and $G$-quasi-invariant $\nu$. Then $G/H$ is a $G$-boundary.
\end{proposition}

Theorem~\ref{thmC:Homogeneous-Poisson-is-Furstenberg} will follow from this proposition: One may use that $S$ is amenable, to invoke the fixed-point property and find some $S$-invariant probability measure $\eta$ on $\partial_F G$. Since $G/S$ is compact, the $G$-orbit of $\eta$ already contains a point mass, which forces $\eta = \delta_x$ for some $x \in \partial_F G$. Then the map $\alpha\colon G/S \twoheadrightarrow \partial_F G, gS \mapsto gx$ is continuous and $G$-equivariant. By the universal property we have $\partial_F G \twoheadrightarrow G/S$ a $G$-equivariant surjection in the other direction. But on every $G$-boundary $X$, the identity map $\operatorname{id}_X$ is the unique continuous, $G$-equivariant self map. So $\alpha$ and $\beta$ are inverses of each other, and $G/S \cong \partial_F G$.

\begin{lemma}\label{lem:quasi-invariant-zero-set}
    Let $\nu$ be $G$-quasi-invariant. If $\nu(A) = 0$ then $\forall x \in G/H \colon m_G(\{g \mid gx \in A\}) = 0$.
\end{lemma}
\begin{proof}
    Let $X = G/H$ and $A \subseteq X$ Borel. If $\nu(A) = 0$, then by Fubini and $G$-quasi-invariance \[
        \int_X m_G(\{g \mid gx \in A\}) \dif \nu(x)
        = \int_G \nu(g^{-1}A) \dif m_G(g)
        = 0.
    \] This shows the claim for some $x \in X$. By $G$-transitivity and $\Delta_G > 0$ we get it for all $x \in X$.
\end{proof}

In particular, Lemma~\ref{lem:quasi-invariant-zero-set} reminds that quasi-invariant $\nu$ have full support: One just notes that for all nonempty, open subsets $U \subseteq G/H$ again $\{g \mid gx \in U\} \subseteq G$ is nonempty, open. And nonempty, open subsets of $G$ have positive Haar measure, so $\nu(U) = 0$ would be a contradiction. So $\supp \nu = G/H$ because its complement is $\bigcup_{U \, \text{open}\colon \nu(U) = 0}U = \emptyset$.

\begin{proof}[Proof of Proposition~\ref{prop:G-boundary-from-poisson}]
    Write $X = G/H$. Being $G$-transitive, $X$ is minimal, so we only have to show strong proximality. In fact, showing some point-measure is in the $G$-orbit of a given probability measure $\eta$ on $X$ will imply that all point masses are contained. So we can assume $\eta$ and $x \in X$ are given, and we have to construct a net $g_i \in G$ such that $g_i \eta \to \delta_x$ in the weak$^*$ sense.
    
    Having that $(X, \nu)$ is a Poisson boundary, the Poisson transform $P_\nu$ is isometric, so for every open neighborhood $U$ of $x$ we have $\operatorname{ess\,sup}_{g \in G} \nu(g^{-1}U) = 1$.
    Since $X$ is metrizable, we can choose a sequence of neighborhoods $x \in U_n \subseteq X$ with radius $1/n$ for some fixed metric, and elements $g_n \in G$ such that $\nu(X \setminus g_n^{-1}U_n) < 1/2^n$. Then by the Borel-Cantelli lemma the following is $\nu$-null: \[
        N = \limsup_{n \to \infty} \,(X \setminus g_n^{-1}U_n) 
        = X \setminus \bigcup_{n \geq 1} \bigcap_{k \geq n} (g_k^{-1}U_k).
    \]
    
    So for $\nu$-a.e. $y \in X$, namely for all $y \in X \setminus N$, we know that $g_n y \to x$. This already means we have $g_n \nu \to \delta_x$ weak$^*$ly. Using the quasi-invariance, this computation transfers to all measures $\eta$ on the homogeneous space, after possibly applying a shift:
    
    Lemma~\ref{lem:quasi-invariant-zero-set} gives $m_G(\{h \in G \mid hy \in N\}) = 0$ for all $y \in X$. By Fubini \[
        \int_G \eta(h^{-1}N) \dif m_G(h)
        = \int_X m_G(\{h \in G \mid hy \in N\}) \dif\eta(y) = 0.
    \] So there exists $h \in G$ with $(h\eta)(N) = 0$. That means $g_n y \to x$ for $(h\eta)$-almost every $y \in X$. This shows $(g_n h) \eta \to \delta_x$ in the weak$^*$ sense using dominated convergence theorem.
\end{proof}

In practice, one may read the combination of Theorems~\ref{thmB:Compact-times-amenable} and~\ref{thmC:Homogeneous-Poisson-is-Furstenberg} as a way to compute Poisson and Furstenberg boundaries for compact-by-amenable groups. However, starting from a factorization $G = KH$ with closed, amenable $H$ one obtains a potentially larger amenable subgroup $S$ with $G = KS$ where $G/S$ identifies the boundaries. The following makes precise how much larger the subgroup $S$ becomes:

\begin{corollary}
    Let $G$ be locally compact, second countable, with compact subgroup $K \leq G$. If $G = KH$ with closed, amenable subgroup $H$, then the Furstenberg boundary of $G$ can be computed: \[
        \partial_F G = G / S, \quad S = \{k \in K \mid \overline{\langle H \cup \{k\} \rangle} \text{ is amenable}\} \cdot H.
    \] Here $\overline{\langle\,\cdot\,\rangle}$ means ``the smallest closed subgroup of $G$ containing $(\,\cdot\,)$''.
\end{corollary}
\begin{proof}
    We can find a cocompact amenable subgroup $S$ such that $\partial_F G \cong G/S$ where $H \subseteq S$, and such that $S$ is maximal in the sense that if any other amenable $A$ contains $H$, then $S$ contains $A$.

    We have $S = (S \cap K) \cdot H$. The inclusion ``$\supseteq$'' holds because $H \subseteq S$. The inclusion ``$\subseteq$'' holds because $G = KH$. Namely, if $S \ni s = kh$ with $k \in K, h \in H$, then already $k = s h^{-1} \in S$.

    We have $S \cap K = \{k \in K \mid \overline{\langle H \cup \{k\}\rangle} \text{ amenable}\}$. The inclusion ``$\subseteq$'' holds because $H \subseteq S$ and because closed subgroups of amenable groups are amenable. The inclusion ``$\supseteq$'' holds because $A = \overline{\langle H \cup \{k\}\rangle}, k \in K$ contains $H$. So if $A$ is amenable then $A \subseteq S$, so $k \in S$.
\end{proof}

\begin{example}
    One may consider $G = \SL_2(\R)$. It contains copies of the compact circle group
    \[
    K = \left\{\m{\cos(\theta) & -\sin(\theta) \\ \sin(\theta) & \cos(\theta)} \mid 0 \leq \theta < 2\pi\right\}
    \]
    and the upper half-plane
    \[
    H = \left\{\m{y & x \\ 0 & 1/y} \mid y > 0\right\},
    \]
    which is a semi-direct product of abelian groups, so solvable and hence amenable. One can write down a unique factorization $G = KH$ with $K \cap H = \{I\}$. 

    Now one computes the Furstenberg boundary as $\partial_F G = G/S$ where $S = \{\pm I\} \cdot H$: In fact, for $k = \pm I$ one checks that $(H \cup - H)$ is closed and still solvable, while for all other $k \neq \pm I$ it is not difficult to write down a formula that expresses a generic element of $G$ in terms of only $k, k^{-1}$ and elements of $H$, but $\langle H, k \rangle = G$ is famously non-amenable.
\end{example}

\section{Examples}\label{sect:examples}

\subsection{Reductive algebraic groups over local fields}\label{sec:red-alg-grps}

Let $\cG$ be a connected reductive algebraic group over an infinite non-Archimedean local field $F$.
We consider the locally compact group $G = \cG(F)$.
Then, by the Bruhat--Tits theory, we obtain a compact open subgroup $K \subset G$ as the stabilizer of an appropriately chosen point in the enlarged Bruhat--Tits building for $\cG$~\cite{MR327923,MR546588}.
In the case of split semisimple groups, we can take $K = \cG(O)$ for the ring of integers $O \subset F$.

By a result of Satake~\cite{Sat63} (see also~\cite[Theorem 4.1]{MR546593}), the spherical Hecke algebra $C_c(K \backslash G / K)$ is commutative.

\begin{proposition}
    Let $\cP \subset \cG$ be the minimal parabolic subgroup defined by a choice of root data for $\cG$. Then the subgroup $S = \cP(F)$ is maximally amenable in $G$.
\end{proposition}

\begin{proof}
    The groups $G$ and $S$ become part of a Tits system, see for example~\cite[paragraph IV.14.15 and Theorem V.21.15]{Bor91}.
    The intermediate subgroups $S \subset P \subset G$ are classified by subsets $I$ of the set of positive simple roots~\cite[Théorème IV.2.5.3]{bourbaki-lie-fr-4-6}.
    If $I$ is not empty, the corresponding parabolic subgroup $P_I$ is indeed nonamenable, since it contains a closed subgroup locally isomorphic to $\SL_2(F)$ in its Levi factor.
\end{proof}

From any initial cocompact amenable subgroup $B$,~\cite[Proposition 4.4]{Fu73} determines that all $G$-boundaries are of the form $G/H$ for some intermediate amenable subgroup $B \subseteq H \subseteq G$.

\begin{corollary}\label{corollary:furstenberg-boundary}
    In the above setting, the Furstenberg boundary of $G$ is $G/S$.
\end{corollary}

\subsection{Type A groups over p-adic numbers}

For the type A groups $G = \SL_n(\Q_p), n \geq 2$ over the $p$-adic numbers $\Q_p$ there are more elementary ways to compute the Furstenberg boundary.

It will be enough to show that $G/B$ is strongly proximal, with $B$ the subgroup of upper triangular matrices of $G$. Because when $G/B$ is a $G$-boundary, from the universal property of $\partial_F G$, we obtain an equivariant surjection $\partial_F G \twoheadrightarrow G/B$. When $\partial_F G \cong G/H$ for $B \subseteq H$ then from the equivariance we get also $H \subseteq g B g^{-1}$ for $g \in G$ with $\phi(eH) = gB$. But the upper triangular subgroup $B$ does not contain any conjugate except itself, so $g^{-1} B g \subseteq H \subseteq B$ implies $g = e$, and thus $H = B$, so $\partial_F G \cong G/B$.

To show the strong proximality of $G \curvearrowright G/B$, one may consider the ``large'' open subset \[
    W = \{ gB \mid \text{all top-left subminors of $g \in G$ are invertible} \} \subseteq G/B.
\] One can show that for every probability measure $\eta$ on $G/B$ one may apply a shift $\widetilde \eta = k\eta, k \in G$ such that $\widetilde \eta$ is supported on $W$. Furthermore, one can show the sequence $a^n \widetilde \eta$ converges weak$^*$ to the point mass $\delta_{eB}$, where $a \in G$ can be any diagonal matrix of the form \[
    a = \operatorname{diag}(p^{\lambda_1}, \dots, p^{\lambda_{n}}),
    \quad \lambda_1 < \dots < \lambda_n, 
    \quad \lambda_1 + \dots + \lambda_n = 0.
\] 

For example, when $n = 2$ the group $G = \SL_2(\Q_p)$ acts on the projective line $G/B \cong \mathbb P^1(\Q_p)$ via the Moebius transform
\[
    \m{a & b \\ c & d} (x : y) = (ax + by : cx + dy).
\]
Here $W = \mathbb P^1(\Q_p) \setminus \{(0:1)\}$, and since $\Q_p$ is uncountable, every probability measure $\eta$ must vanish on some $(x:y)$, so one can translate its support into $W$. Then positive powers of
\[
a = \m{p^{-1} & 0 \\ 0 & p^{+1}}
\]
contract $W$ to $\{(1 : 0)\}$.

\printbibliography

\end{document}